\documentclass[11pt]{article}
\title{\vskip-1.5em Approximate amenability of Schatten classes, Lipschitz algebras and second duals of Fourier algebras}
\author{Y. Choi, F. Ghahramani\thanks{Research supported by NSERC grant 36640-07}}
\date{11th June 2009}
\date{1st October 2009}
\usepackage{mathpazo}

\RequirePackage{amsmath,amsfonts,amsthm}

\newcommand{\appctr}{approx\-imately con\-tract\-ible}
\newcommand{\appam}{approx\-imately amen\-able}
\newcommand{\appctry}{approx\-imate con\-tract\-ib\-ility}
\newcommand{\appamy}{approx\-imate amen\-ab\-ility}

\newcommand{\dt}[1]{{\it #1}\/}
\newcommand{\st}{\,:\,}
\newcommand{\fu}[1]{{#1}^\#}
\newcommand{\defeq}{:=}

\newcommand{\iso}{\cong}
\newcommand{\boxprod}{\mathop{\Box}\nolimits}  

\newcommand{\wstar}{\mathop{\ensuremath{\text{weak}^*}}\nolimits}
\newcommand{\tp}{\mathop{\otimes}\nolimits}
\newcommand{\ptp}{\mathop{\widehat{\otimes}}\nolimits}
\newcommand{\lp}[1]{\ell^{#1}}

\newcommand{\id}[1][]{{\mathsf 1}_{#1}}
\newcommand{\abs}[1]{{\left\vert#1\right\vert}}
\newcommand{\norm}[1]{{\Vert#1\Vert}}
\newcommand{\Norm}[1]{{\left\Vert#1\right\Vert}}

\newcommand{\al}{\alpha}
\newcommand{\veps}{\varepsilon}

\newcommand{\cA}{{\mathcal A}}
\newcommand{\cB}{{\mathcal B}}
\newcommand{\cF}{{\mathcal F}}
\newcommand{\cS}{{\mathcal S}}

\newcommand{\Nat}{{\mathbb N}}
\newcommand{\Nmin}{{\mathbb N}_{\min}}

\newcommand{\Cplx}{{\mathbb C}}
\newcommand{\Real}{{\mathbb R}}

\newcommand{\Rad}{\operatorname{Rad}}
\newcommand{\supp}{\operatorname{supp}}

\newcounter{pulse}
\numberwithin{pulse}{section}
\newcommand{\theadft}[1]{\sc #1}
\theoremstyle{plain}
\newtheorem{thm}[pulse]{\theadft{Theorem}}
\newtheorem{prop}[pulse]{\theadft{Proposition}}
\newtheorem{lem}[pulse]{\theadft{Lemma}}

\theoremstyle{definition}
\newtheorem{defn}[pulse]{\theadft{Definition}}
\newtheorem{rem}[pulse]{\theadft{Remark}}
\newtheorem{rems}[pulse]{\theadft{Remarks}}
\newtheorem{eg}[pulse]{\theadft{Example}}

\usepackage[usenames]{color}
\usepackage{hyperref}

\numberwithin{equation}{section}

\renewcommand{\dt}[1]{\textcolor{Bittersweet}{\sf {#1}}}  
\newcommand{\mnorm}[1]{\norm{#1}_{\rm mul}}  

\newcommand{\dlVP}{de la Vall\'ee Poussin}
\newcommand{\fourier}{A(\widehat{\Real^d})}

\newcommand{\snorm}[1]{p_\al(#1)}
\newcommand{\lip}{\operatorname{lip}}
\newcommand{\Lip}{\operatorname{Lip}}

\newcommand{\FA}{\operatorname{A}}
\newcommand{\FS}{\operatorname{B}}
\newcommand{\VN}{\operatorname{VN}}
\newcommand{\Cstar}{\operatorname{C}^*}
\newcommand{\UCBhat}[1]{\operatorname{UCB}(\widehat{#1})}
\newenvironment{YCnum}{%
\begin{enumerate}

}{\end{enumerate}\ignorespacesafterend}

\begin{document}
\maketitle

\begin{abstract}
Amenability of any of the algebras described in the title is known to force them to be finite-dimensional. The analogous problems for \emph{approximate} amenability have been open for some years now.
In this article we give a complete solution for the first two classes, using a new criterion for showing that certain Banach algebras without bounded approximate identities cannot be approximately amenable. The method also provides a unified approach to existing non-approximate amenability results, and is applied to the study of certain commutative Segal algebras.

Using different techniques, we prove that \emph{bounded} approximate amenability of the second dual of a Fourier algebra implies that it is finite-dimensional. Some other results for related algebras are obtained.

\end{abstract}

\begin{section}{Introduction}
The concepts of \appamy, \appctry\ and essential amen\-ab\-ility, and some other related concepts, were introduced and studied in \cite{GhL_genam1} and further developed in \cite{GhSt_AG,GhZ_pseudo,DLZ_St06,GhLZ_genam2,CGZ_JFA09}. 
The studies include determining when various classes of Banach algebras are, or are not, \appam.
Attempting to prove that a given Banach algebra is \emph{not} \appam\ can be much more difficult than for `classical' amenability. 

It was shown in \cite{DLZ_St06} that for $1\leq p<\infty$ the algebras $\ell^p$, equipped with pointwise multiplication, are not \appam;
this has motivated
the question of whether the Schatten classes $\cS^p(H)$ are approximately amenable.
 In \cite[\S7]{GhLZ_genam2}, partial results were also obtained about \appamy, or otherwise, of the little Lipschitz algebras $\lip_\alpha(X)$. 

It remains an open question, to the authors' knowledge, if an \appam\ Banach algebra $\cA$ must necessarily have a bounded approximate identity.
The authors have recently proved, in collaboration with Y. Zhang, that the {\it a priori}\/ stronger notion of bounded \appctry\ forces the existence of a bounded approximate identity; see \cite[\S3]{CGZ_JFA09}. 
The main result of Section~\ref{s:criterion}
gives a criterion for
proving that certain algebras without bounded approx\-imate identities
cannot be \appam. 

For example, we use our criterion to show that neither Schatten classes on infinite-dimensional Hilbert spaces, nor Lipschitz algebras on infinite compact metric spaces, are approximately amenable. Another question left open from previous studies is whether every proper Segal subalgebra of the group algebra of a locally compact group fails to be \appam. We prove that this is the case for the group~$\Real^n$. We correct an error in the literature on \appamy\ of the semigroup algebras of Brandt semigroups, and also use our criterion to provide new proofs for several results obtained by other authors.  

Another open question has been whether every \appam\ Banach algebra whose underlying space is reflexive, as a Banach space, is finite-dimensional.  We show this is the case, under the additional hypotheses of commutativity and \emph{bounded} \appamy.
Finally we show that the second dual algebra of the Fourier algebra of a locally compact group is boundedly \appam, if and only if the group is finite. This extends the result in the amenable case, which was proved by A.~T. Lau and R.~J. Loy \cite{LauLoy_JFA97} and independently by E.~Granirer~\cite{Gran_JAust97}.

\paragraph\ Now we recall some basic definitions. A derivation $D$ from a Banach algebra $\cA$ into a Banach $\cA$-bimodule $X$ is \dt{approximately inner}, if there exists a net $(\operatorname{ad}_{x_i})$ of inner derivations from $\cA$ into $X$ such that $D$ is the strong-operator-topology limit of $(\operatorname{ad}_{x_i})$\/.
The Banach algebra $\cA$ is \dt{\appam} if every continuous derivation $D$ from $\cA$ into any dual Banach $\cA$-bimodule $X^*$ is approximately inner, for all Banach $\cA$-bimodules $X$. If in the above definition the nets of inner derivations  can always be taken to be bounded, then $\cA$ is \dt{boundedly \appam}.

Approximate amenability and bounded \appamy\ can be characterized by the existence of so-called ``approximate diagonals'' (\cite{GhL_genam1,GhLZ_genam2}), and this will be used several times in the proofs of our main results.   
\end{section}

\begin{section}{A criterion for ruling out \appamy}\label{s:criterion}
It was shown in \cite[\S4]{CGZ_JFA09} that the Fourier algebra of the free group on two generators is not (operator) \appam. Motivated by some of the ideas in that argument, we can now give a criterion which shows that numerous examples of Banach algebras without bounded approximate identities fail to be \appam.

In particular, we have the following result, which resolves a question that had been open since~\cite{GhL_genam1}. 
(Recall from \cite[Corollary~7.1]{GhL_genam1} that $\cS_p(H)$ is \emph{essentially amenable}.)

\begin{thm}\label{t:Schatten}
Let $1\leq p <\infty$ and let 
$\cS_p(H)$ denote the ideal of Schatten class operators on an infinite-dimensional Hilbert space~$H$. Then $\cS_p(H)$ is not \appam.
\end{thm}

We shall deduce this from a much more general result, whose statement requires some technical definitions.
Given a Banach algebra $\cA$\/,
we define the \dt{multiplier norm on $\cA$} by $\mnorm{a}\defeq \max(\norm{\lambda_a},\norm{\rho_a})$\/, where $\lambda_a: \cA\to \cA$\/, $x\mapsto ax$ is left multiplication by $a$ and $\rho_a: \cA\to \cA$\/, $x\mapsto xa$ is right multiplication by~$a$\/.

\begin{defn}\label{dfn:separated-lumps}
A \dt{separated, unbounded, multiplier-bounded configuration} (or \dt{SUM} configuration for short) in $\cA$ consists of two sequences $(u_n), (p_n)\subseteq \cA$
which satisfy the following properties.
\begin{itemize}
\item[] (Separated) $u_np_n=p_nu_n=u_n$ for all $n$\/; and $u_jp_k=p_ku_j=0$\/, whenever $j\neq k$\/.
\item[] (Unbounded) $\norm{u_n}\to\infty$ as $n\to\infty$\/.
\item[] (Multiplier-bounded) The sequences $(\mnorm{u_n})$ and $(\mnorm{p_n})$ are bounded.
\end{itemize}
\end{defn}

\begin{rems}\
\begin{YCnum}
\item In the definition of a SUM configuration, by passing to a subsequence if necessary, we may assume that the sequence $(\norm{u_n})_{n\geq 1}$ grows faster than any prescribed sequence of positive reals.
\item If $\cA$ is a Banach algebra with a bounded approximate identity, then $\norm{\cdot}$ and $\mnorm{\cdot}$ are equivalent norms, by Cohen's factorization theorem, so that $\cA$ contains no SUM configurations.
\item In several, though not all, of the applications below, one can choose the $(u_n)$ to be pairwise orthogonal idempotents, and take $p_n=u_n$\/.
\end{YCnum}
\end{rems}

The definition of a SUM configuration may seem overly technical. It is therefore convenient to have a simpler criterion that ensures the existence of such sequences.

\begin{lem}\label{l:absorbing-unbounded}
Let $\cA$ be a Banach algebra. Suppose that there exists an unbounded but multiplier-bounded sequence $(E_n)_{n\geq 1}\subseteq \cA$\/,
such
that $E_nE_{n+1}=E_n=E_{n+1}E_n$ for all~$n$\/. Then $\cA$ contains a SUM configuration.
\end{lem}

\begin{proof}[Proof of Lemma~\ref{l:absorbing-unbounded}]
Note that, for $1\leq j<k$\/, we have $E_jE_k=E_j=E_kE_j$\/. Moreover, by passing to a subsequence if necessary, we may assume without loss of generality that $\norm{E_{n+1}} > \norm{E_n}$ for all $n$ and that $\norm{E_{n+1}} - \norm{E_n} \to \infty$ as $n\to\infty$\/.

We now define
\[ u_n\defeq E_{3n}-E_{3n-1} \quad,\quad p_n\defeq E_{3n+1} - E_{3n-2} \quad\text{ for $n\geq 1$\/.} \]
Clearly both $(\mnorm{u_n})$ and $(\mnorm{p_n})$ are uniformly bounded sequences, since the sequence $(E_n)$ is multiplier-bounded. Moreover, our construction ensures that
\[ \norm{u_n} \geq \norm{E_{3n}}-\norm{E_{3n-1}} \to\infty \,.\]
It remains to verify the `separatedness' conditions.
Since the set $\{E_n\;:\;n\in\Nat\}$ is commutative, $u_j$ and $p_k$ commute for every $j,k$\/. 
 We have
\[ u_np_n = p_nu_n = (E_{3n+1} - E_{3n-2})(E_{3n}-E_{3n-1})
 = E_{3n} - E_{3n-1}  = u_n \quad\text{ for all $n$\/.} \]
If $j < k$\/, then, since $3k-2>3j$\/,
\[ u_jp_k = p_k u_j = (E_{3k+1}-E_{3k-2})(E_{3j}-E_{3j-1}) = 0\,;\]
and similarly, since $3k-1>3j+1$\/,
\[ u_kp_j = p_ju_k = (E_{3j+1} - E_{3j-2})(E_{3k}-E_{3k-1}) = 0 \,.\]
Thus $(u_n)$ and $(p_n)$ form a SUM configuration, as required.
\end{proof}

The point of our definition is the following result, whose proof should be compared with the main arguments in \cite[\S4]{CGZ_JFA09}. 
\begin{thm}\label{t:better-nonAA}
Let $\cA$ be a Banach algebra which contains a SUM configuration.
Then $\cA$ is not \appam.
\end{thm}

\begin{proof}[Proof of Theorem \ref{t:better-nonAA}]
Let $(u_n)$, $(p_n)$ be a SUM configuration in $\cA$, and put
\[ C:= \max\{\sup\nolimits_n \mnorm{u_n},\sup\nolimits_n \mnorm{p_n}\}\,.\]
Elementary calculations using
the definition of
the projective tensor norm show that for any $w\in \cA\ptp \cA$\/,
\begin{equation}\label{eq:mbdd}
\max( \norm{u_n\cdot w }, \norm{w\cdot u_n},  \norm{p_n\cdot w }, \norm{w\cdot p_n} ) \leq C\norm{w} \quad\text{ for all $n$\/.} \tag{$*$}
\end{equation}

By passing to a subsequence if necessary, we can without loss of generality assume that in addition to satisfying \eqref{eq:mbdd}, our sequence $(u_n)$ satisfies
\begin{equation}\label{eq:growing-sufficiently-fast}
\norm{u_n} \geq (n+1)^4 \quad\text{ for all $n$\/.}  \tag{$\dagger$}
\end{equation}
(The exact rate of growth is not important, but we need it to be moderately fast.)

Note that the `separatedness conditions' imply that the $(u_n)$ are pairwise orthogonal,
 since for $j\neq k$ we have $u_ju_k = u_jp_ju_k=0$\/.
Define $a,b\in \cA$ by
\[ a= \sum_{j\geq 1} (2j-1)^{-2} \norm{u_{2j-1}}^{-1} u_{2j-1} 
	\quad\text{ and }\quad
 b= \sum_{k\geq 1} (2k)^{-2} \norm{u_{2k}}^{-1} u_{2k}\,, \]
and observe that for each $n\in\Nat$\/:
\begin{equation}\label{eq:orthog}
\begin{aligned}
 p_{2n-1}a & = (2n-1)^{-2} \norm{u_{2n-1}}^{-1} u_{2n-1} \,;
 & au_{2n} & = 0 \,; \\
 bp_{2n} &= (2n)^{-2} \norm{u_{2n}}^{-1} u_{2n} \,;
 & u_{2n-1}b & = 0\,.
\end{aligned}
\tag{$**$}
\end{equation}

Our proof now proceeds by contradiction. Suppose that $\cA$ is \appam, and let $\veps>0$\/. Then, since \appamy\ and \appctry\ are the same,
there exist $M\in \cA\ptp \cA$ and $F,G\in \cA$ such that:
\begin{subequations}
\begin{align}
\label{eq:main_a}
\norm{ a\cdot M - M\cdot a - a \tp G + F \tp a} & < \veps \;; \\
\label{eq:main_b}
\norm{ b\cdot M - M\cdot b - b \tp G + F \tp b} & < \veps \;;
\end{align}
\end{subequations}
\begin{equation}\label{eq:near_1}
\norm{a-aF} < \veps\,.
\end{equation}

From \eqref{eq:main_a} and the observations \eqref{eq:mbdd} and \eqref{eq:orthog}, we see that for any $j,k\in\Nat$,
\[ \begin{aligned}
C^2\veps
	 & > \norm{p_{2j-1}\cdot( a\cdot M - M\cdot a - a \tp G + F \tp a  )\cdot u_{2k}} \\
	 & = \norm{p_{2j-1} a\cdot M\cdot u_{2k} - p_{2j-1}a \tp Gu_{2k} } \\
	 & = \norm{ (2j-1)^{-2} \norm{u_{2j-1}}^{-1}
	(u_{2j-1}\cdot M\cdot u_{2k} - u_{2j-1}\tp Gu_{2k}) }\/,
\end{aligned} \]
so that
\begin{subequations}
\begin{equation}\label{eq:next_a}
\frac{\norm{ u_{2j-1}\cdot M\cdot u_{2k} - u_{2j-1}\tp Gu_{2k} } }{\norm{u_{2j-1}}\,\norm{u_{2k}}} \leq C^2\veps \frac{(2j-1)^2}{\norm{u_{2k}} } \;.
\end{equation}
A similar argument, using \eqref{eq:main_b} instead of \eqref{eq:main_a}, yields
\begin{equation}\label{eq:next_b}
\frac{\norm{ - u_{2j-1}\cdot M\cdot u_{2k} + u_{2j-1}F\tp u_{2k}  } }{\norm{u_{2j-1}}\,\norm{u_{2k}}} \leq C^2\veps \frac{(2k)^2}{\norm{u_{2j-1}} } \;;
\end{equation}
\end{subequations}
combining \eqref{eq:next_a} and \eqref{eq:next_b} using the triangle inequality then gives
\[ 
\begin{aligned}
   C^2\veps \left( \frac{(2j-1)^2}{\norm{u_{2k}}} + \frac{(2k)^2}{\norm{u_{2j-1}}} \right)
  & \geq  \Norm{ \frac{u_{2j-1}F\tp u_{2k}}{\norm{u_{2j-1}}\,\norm{u_{2k}} }
 - \frac{u_{2j-1}\tp Gu_{2k}}{\norm{u_{2j-1}}\,\norm{u_{2k}}} } \\
  & \geq  \abs{ \frac{\norm{u_{2j-1}F}}{\norm{u_{2j-1}}} - \frac{\norm{Gu_{2k}}}{\norm{u_{2k}}} }
\end{aligned}
\]
(Here we used the inequality $\norm{x-y} \geq \abs{ \norm{x}-\norm{y} }$, valid in any normed space, together with the fact that the projective tensor norm is a cross-norm). Using the growth condition \eqref{eq:growing-sufficiently-fast}, we therefore have
\begin{equation}\label{eq:decoupled}
\abs{ \frac{\norm{u_{2j-1}F}}{\norm{u_{2j-1}}} - \frac{\norm{Gu_{2k}}}{\norm{u_{2k}}} }
\leq   C^2\veps \left( \frac{(2j-1)^2}{(2k+1)^4} + \frac{(2k)^2}{(2j)^4} \right).
\end{equation}
Taking $j=k$ in \eqref{eq:decoupled} gives
\[ \abs{ \frac{\norm{u_{2k-1}F}}{\norm{u_{2k-1}}} - \frac{\norm{Gu_{2k}}}{\norm{u_{2k}}} }
    \leq C^2\veps \left( \frac{(2k-1)^2}{(2k+1)^4} + \frac{(2k)^2}{(2k)^4}\right) 
    \leq \frac{2C^2}{(2k)^2}\veps\,,
\]
while taking $j=k+1$ in \eqref{eq:decoupled} gives
\[
 \abs{ \frac{\norm{u_{2k+1}F}}{\norm{u_{2k+1}}} - \frac{\norm{Gu_{2k}}}{\norm{u_{2k}}} }
    \leq C^2\veps \left( \frac{(2k+1)^2}{(2k+1)^4} + \frac{(2k)^2}{(2k+2)^4}\right) 
    \leq \frac{2C^2}{(2k+1)^2}\veps\,.
\]
Combining the preceding two inequalities, we see that, for every $k\geq 1$\/,
\[ \abs{\frac{\norm{u_{2k-1}F}}{\norm{u_{2k-1}}} - \frac{\norm{u_{2k+1}F}}{\norm{u_{2k+1}}} } \leq 2C^2\veps \left( \frac{1}{(2k)^2} + \frac{1}{(2k+1)^2} \right)\;. \]
Then, since $\lim_{m\to\infty} \norm{u_{2m+1}}^{-1} \norm{u_{2m+1}F}=0$\/, we deduce that
\begin{equation}\label{eq:telescope}
\begin{aligned}
\frac{\norm{u_1F}}{\norm{F}}
  & = \lim_{m\to\infty} \abs{ \frac{\norm{u_1F}}{\norm{u_1}} - \frac{ \norm{u_{2m+1}F} }{\norm{u_{2m+1}}} } \\ 
  & = \lim_{m\to\infty} \abs{
	\sum_{k=1}^m \frac{\norm{u_{2k-1}F}}{\norm{u_{2k-1}}} - \frac{\norm{u_{2k+1}F}}{\norm{u_{2k+1}}} 
} \\
  & \leq 2C^2\veps\, \sum_{k=1}^\infty  \frac{1}{(2k)^2} + \frac{1}{(2k+1)^2} 
  & \leq 2C^2\veps\,.
  \end{aligned}
\end{equation}

However: going back to the `approximate identity' formula \eqref{eq:near_1}, we also have
\[
 C\veps \geq \norm{ u_1(a-aF)} = \Norm{ \frac{u_1}{\norm{u_1}} - \frac{u_1F}{\norm{u_1}} } \geq \abs{ 1 - \frac{\norm{u_1F}}{\norm{u_1}} } \;; \]
combining this with \eqref{eq:telescope},
 we conclude that
\[ 1\leq C\veps + 2C^2\veps\,. \]
Since $\veps$ was arbitrarily small this yields a contradiction, and therefore $\cA$ is not \appam.
\end{proof}

\begin{proof}[Proof of Theorem~\ref{t:Schatten}]
We can now complete the proof that $\cS_p(H)$ is not \appam. While we could use Lemma~\ref{l:absorbing-unbounded}, it is more instructive to 
construct an explicit SUM configuration as follows.
Fix an
infinite orthonormal set
 $(\delta_n)_{n\geq 1}$ and take $P_1$ to be the orthogonal projection onto the linear span of $\delta_1$\/; take $P_2$ to be the orthogonal projection onto the linear span of $\{\delta_2, \delta_3\}$\/; and so on. Since $\norm{T}_\cA =\left(\sum_{k=1}^\infty s_k(T)^p\right)^{1/p}$\/, where
\[  s_k(T)= \inf \{ \norm{T-R}_{\rm op} \st \operatorname{rank}(R) < k \}\,, \]
it is easily checked that $\mnorm{P_n}=1$\/. On the other hand, since $\norm{T}_\cA= \left(\mathop{\rm trace} \abs{T}^p\right)^{1/p}$\/, we have $\norm{P_n}_\cA = n^{1/p}$ for all~$n$\/. As the $(P_n)$ are pairwise orthogonal idempotents, taking $p_n=u_n=P_n$ yields a SUM configuration. We finish by applying Theorem~\ref{t:better-nonAA}.
\end{proof}

\begin{rem}
The same argument ought to work for other operator ideals on a large class of Banach spaces, although some hypothesis on the space is necessary in order to produce sufficiently many projections.
\end{rem}

In the remainder of this section, we shall give some further examples of how our criterion may be applied.

\begin{eg}
Let $(B_n)$ be a sequence of Banach algebras, such that for each $n$ the algebra $B_n$ has an identity element $\id[n]$\/, and suppose that $\norm{\id[n]}_{B_n}\to\infty$ as $n\to\infty$\/. Let $\cA$ be the completion of the algebraic direct sum $\bigoplus_n B_n$ in either the $c_0$ or $\ell^p$ norms for $1\leq p <\infty$\/.
When regarded as an element of $\cA$, each $\id[n]$ is an idempotent with norm~$\norm{\id[n]}_{B_n}$ and multiplier norm~$1$\/. Hence on taking $u_n=p_n=\id[n]$ we have produced a SUM configuration in $\cA$\/, and so by Theorem \ref{t:better-nonAA} $\cA$ is not \appam.

In particular, we see that the $c_0$-sum of a sequence $(B_n)$ of finite-dimensional amenable Banach algebras need not be approximately amenable (just take each $B_n$ to be $\Cplx^n$ with pointwise multiplication and the $\ell^1$-norm). This should be compared with \cite[Example 6.1]{GhL_genam1}.
\end{eg}

\begin{eg}[Correction of \cite{LashSam_SF05}]
Let $S$ be a Brandt semigroup over a group $G$ with \emph{infinite} index set $I$\/, and consider the convolution algebra $\ell^1(S)$\/. It is claimed in \cite[Theorem 1.8]{LashSam_SF05} that if $G$ is amenable then $\ell^1(S)$ will be \appam: however, the proof contains an error. (Briefly: in the notation of \cite{LashSam_SF05}, the algebras $\ell^1(S_F)$ are in general \emph{not} unital, since the norm of the identity element in $\ell^1(S_F)$ grows with $\abs{F}$; hence the appeal to results of \cite{GhL_genam1} is invalid.)

Indeed, $\ell^1(S)$ is \emph{not} \appam. For, let $\cA$ be the Banach algebra which consists of all $\Cplx$-valued $I\times I$ matrices $(a_{ij})_{i,j\in I}$ such that $\sum_{i,j\in I} \abs{a_{ij}}< \infty$\/, equipped with the usual matrix multiplication as its product. As observed in \cite{DuncNam} (see the remarks before Theorem 8), the algebra $\ell^1(S)$ has a 1-dimensional ideal $\Cplx\theta$ (where $\theta$ is a formal zero element of $S$), and the quotient algebra $\ell^1(S)/\Cplx\theta$ is isomorphic to $\cA\ptp \ell^1(G)$\/. If we let $\veps: \ell^1(G)\to \Cplx$ be the augmentation character, then the composite map
\[ \ell^1(S) \to \ell^1(S)/\Cplx\theta \iso \cA\ptp\ell^1(G)  \xrightarrow{{\sf id}\tp \veps} \cA \]
is a quotient homomorphism of Banach algebras.

Since \appamy\ is inherited by quotients, if $\ell^1(S)$ were \appam, then $\cA$ would be too. This is not the case, because it contains SUM configurations.
For instance, let $(F_n)_{n\geq 1}$ be a sequence of pairwise disjoint, finite subsets of~$I$ such that $\abs{F_n}\to\infty $ as $n\to\infty$\/. Define $u_n$ to be the $I\times I$ matrix with entries given by
\[ (u_n)_{ij} = \left\{ \begin{aligned}
	1  & \quad\text{ if $i=j\in F_n$} \\
	0  & \quad\text{ otherwise.}
	\end{aligned} \right. \]
Then: (i) $\norm{u_n}_\cA = \abs{F_n}$ for all $n$\/; (ii) the $(u_n)$ are pairwise orthogonal idempotents; and (iii) $\mnorm{u_n}=1$ for all $n$\/. Hence, taking $p_n=u_n$ yields a SUM configuration in~$\cA$ as required.
\end{eg}

\begin{rem}
After the present work was done, we were informed that a joint paper of Maysami-Sadr and Pourabbas also contains a proof that the preceding example is not \appam: see~\cite[Theorem~4.5]{Sadr-Pour}. This is deduced from a more general result, \cite[Theorem 3.4]{Sadr-Pour}, which we suspect might also be provable using our techniques; however, we have not attempted to do this.
\end{rem}

\paragraph\
The final two examples of this section are already known not to be \appam. Nevertheless, Theorem~\ref{t:better-nonAA} not only provides alternative proofs, but does so in a unified treatment.
In both cases it is convenient to use Lemma~\ref{l:absorbing-unbounded}.

\begin{eg}[A result of Dales and Loy]
Let $\Nat_{\min}$ denote the semilattice whose underlying set is $\Nat$ and where the `product' of two elements is defined to be their minimum. Any function $\omega:\Nat \to [1,\infty)$ satisfying $\omega(1)=1$ defines a weight function on $\Nat_{\min}$\/, and the weighted convolution algebra $\lp{1}(\Nat_{\min},\omega)$ is then isomorphic (via the Gelfand transform) to the \dt{Feinstein algebra} $A_\omega$\/.

If $\liminf_n \omega(n) < \infty$ then $\lp{1}(\Nat_{\min},\omega)$ has a bounded approximate identity, and hence  by results of \cite{DLZ_St06} (or a direct construction) it is boundedly \appctr. On the other hand, in the case where $\omega(n)\to\infty$ as $n\to\infty$\/, H.~G.\ Dales and R.~J.\ Loy have shown that $\lp{1}(\Nat_{\min},\omega)$ is not even \appam\ (private communication).

The argument of Dales and Loy proceeds by direct calculations, patterned on the arguments of~\cite{DLZ_St06}.
We can now give an alternative proof.
Let $e_j$ denote the point mass concentrated at $j$\/, regarded as an element of $\lp{1}(\Nat_{\min},\omega)$\/.
Since $\lim_j \omega(j)=+\infty$\/, we can inductively construct a strictly increasing subsequence $(j_n)_{n\geq 1} \subset \Nat$ such that
\begin{equation}\label{eq:uphill}
 \text{ for each $n$\/ and every $k\geq j_n$, $\omega(k)\geq\omega(j_n)$.}
\tag{$\spadesuit$}
\end{equation}
Put $E_n = e_{j_n}$\/; then clearly the sequence $(E_n)$ is unbounded in $\lp{1}(\Nat_{\min},\omega)$\/, and satisfies $E_nE_{n+1}=E_n$ for all~$n$\/. Moreover,
 given $a= \sum_i a_i e_i\in \ell^1(\Nat_{\min},\omega)$\/, we have
\[ \begin{aligned}
\norm{E_na}_{1,\omega}
  = \Norm{ \sum_{i=1}^{j_n-1} a_ie_i + \left(\sum_{k \geq j_n} a_k\right)e_{j_n} } 
 & = \sum_{i=1}^{j_n-1} \abs{a_i}\omega(i) + \abs{\sum_{k \geq j_n} a_k}\omega(j_n) \\
 & \leq \sum_{i=1}^{j_n-1} \abs{a_i}\omega(i) + \sum_{k \geq j_n}\abs{a_k}\omega(j_n) \\
 & \leq \sum_{i=1}^{j_n-1} \abs{a_i}\omega(i) + \sum_{k \geq j_n}\abs{a_k}\omega(k)  
 & = \norm{a}_{1,\omega}\;,
\end{aligned} \]
where we used the condition \eqref{eq:uphill} at the penultimate step. Thus, $\mnorm{E_n}=1$ for all~$n$\/; and by applying Lemma~\ref{l:absorbing-unbounded} and Theorem~\ref{t:better-nonAA}, we conclude that $\lp{1}(\Nat_{\min},\omega)$ is not \appam.
\end{eg}

\begin{eg}[A result of P. Lawson]
Let $X$ be a Banach space with a normalized unconditional basis $(\delta_n)$\/. It is standard (see, for instance, \cite[Propositions 1.c.6 and 1.c.7]{LindTzaf_v1}) that $X$ may be equipped with an equivalent norm $\norm{\cdot}$ that has the following two properties:
\begin{itemize}
\item $\norm{\delta_n}=1$ for all $n$\/;
\item for every sequence $(\al_n)\subset \Cplx$ such that the series $\sum_n \al_n \delta_n$ converges in $X$\/, and every sequence $(\gamma_n)\subset\Cplx$ such that $\abs{\gamma_n}\leq \abs{\al_n}$ for all $n$\/, the series $\sum_n\gamma_n \delta_n$ converges in $X$ \emph{and}
\begin{equation}\label{eq:shrinking}
 \norm{\sum_n \gamma_n \delta_n } \leq \norm{\sum_n \al_n \delta_n } \;.
\end{equation}
\end{itemize}
Equipped with $\norm{\cdot}$, $X$ thus
becomes a commutative Banach algebra with respect to pointwise product of sequences. Note that the multiplier norm of $\sum_n \al_n \delta_n$ is then at most $\sup_n \abs{\al_n}$\/.

If $X=c_0$\/, it is amenable and hence {\it a~fortiori}\/ \appam. In contrast, Lawson has shown \cite{PDL_PhD} that if $X$
is not isomorphic (as a Banach algebra) to~$c_0$\/, then it is not even \appam.
His proof, which generalizes arguments from~\cite{DLZ_St06}, starts by making use of the following observation.

\begin{lem}\label{l:not-c_0}
Suppose that $X$ is not algebra-isomorphic to $c_0$\/. Then the function $n\mapsto \norm{ \sum_{i=1}^n \delta_i}$\/, which by \eqref{eq:shrinking} is non-decreasing, must be unbounded.
\end{lem}
For sake of completeness, we sketch a proof.

\begin{proof}[Sketch of proof]
Suppose that $\sup_n \norm{\sum_{i=1}^n \delta_i} = K < \infty$\/. We derive a contradiction as follows: given $a\in c_{00}$\/, let $N=\max\{i\st a_i\neq 0\}$\/, and let $M$ be such that $\abs{a_M}=\norm{a}_\infty$\/. Then \eqref{eq:shrinking} implies that
\[ \norm{a}_\infty = \norm{a_M\delta_M}\leq \norm{\sum_i a_i \delta_i} \leq \norm{a}_\infty \norm{\sum_{i =1}^N \delta_i } \leq K \norm{a}_\infty \,.\]
It follows that the map $a\mapsto \sum_i a_i \delta_i$ extends continuously to a Banach space isomorphism $\theta: c_0 \to X$\/, contradicting the assumption on~$X$\/.
\end{proof}

In the present context, note that Lemma~\ref{l:not-c_0} immediately furnishes a sequence satisfying the condition of Lemma~\ref{l:absorbing-unbounded} (just take partial sums of the form $\sum_{i=1}^n \delta_i$\/). Hence by invoking Lemma~\ref{l:absorbing-unbounded} and Theorem~\ref{t:better-nonAA}, we obtain an alternative proof of Lawson's result. (For special cases such as $X=\ell^p$\/, for $1\leq p <\infty$\/, it is simpler to construct an explicit  SUM configuration rather than go through Lemma~\ref{l:absorbing-unbounded}.)
\end{eg}

\end{section}

\begin{section}{Further applications of our criterion to certain function algebras}\label{s:application}

We now give some other examples covered by Theorem~\ref{t:better-nonAA}, which to the authors' knowledge are new.

Let $G$ be a locally compact group and $L^1(G)$ its group algebra. It has been an open question whether every proper Segal subalgebra of $L^1(G)$\/, as defined in 
\cite[Ch.~6]{Reit_Steg},
fails to be \appam. Recently, H.~G. Dales and R.~J. Loy have provided various examples of Segal algebras on abelian groups, that fail to be \appam\ (private communication). It was already known \cite{CGZ_JFA09} that no proper Segal subalgebra of $L^1(G)$ can be \emph{boundedly} \appam.

\subsection*{Segal subalgebras of $L^1(\Real^d)$}

\begin{thm}\label{t:segal-Rd}
Let $d\in\Nat$\/, and let $\cA$ be a proper Segal subalgebra of
$L^1(\Real^d)$\/. Then $\cA$ is not \appam.
\end{thm}

\begin{proof}
It is slightly more convenient for us to take Fourier transforms and work in the Fourier algebra $\fourier$\/. Let $\cF: L^1(\Real^d)\to \fourier$ be the Fourier transform, normalized so as to be an isometric algebra isomorphism. Then $\cF$ maps the Segal subalgebra $\cA$ onto an `abstract Segal subalgebra' of $\fourier$\/, and we shall denote this image algebra by $\cF\cA$ (equipped with the norm induced from $\cA$ via $\cF$).
 The Segal algebra condition implies that there exists $C>0$ such that
\begin{equation}\label{eq:Segal-condition} 
\norm{gh}_{\fourier}\leq \norm{gh}_{\cF\cA} \leq C \norm{g}_{\fourier} \norm{h}_{\cF\cA} \qquad\text{ for all $g\in \fourier$, $h\in \cF\cA$\/.}
\end{equation}

We shall make crucial use of the \dt{\dlVP\ kernel} for $L^1(\Real^d)$, see \cite[VI.1.13]{KatzHA_3rd}\/. Given $r>0$\/, let $b_r:\Real\to [0,1]$ be the function defined by 
\[ b_r(x) = \left\{ \begin{aligned}
 1 & \quad\text{ if $0\leq \abs{x}\leq r$\/;} \\
 2-r^{-1}x & \quad\text{ if $r\leq \abs{x}\leq 2r$\/;}\\
 0 & \quad\text{ if $\abs{x}\geq 2r$\/.}
\end{aligned}\right. \]
We then define $V_r \in C_c(\widehat{\Real^d})\subset \fourier$ by
\[ V_r(\mathbf{x}) = b_r(x_1)\dotsb b_r(x_d) \qquad(\mathbf{x}=(x_1,\ldots,x_d)\in\Real^d).\]
We note that $V_r$ is equal to $1$ on the `cube' $[-r,r]^d$\/, and vanishes outside the cube $[-2r,2r]^d$\/.

It is known that $\norm{\cF^{-1}b_r}_{L^1(\Real^d)}\leq 3$ for all $r$ (this follows by rewriting $\cF^{-1}b_r$ as a combination of F\'ejer kernels), and therefore $\norm{V_r}_{\fourier}\leq \norm{b_r}^d\leq 3^d$\/, for all~$r$\/.
Therefore \eqref{eq:Segal-condition} implies that
\begin{equation}\label{eq:dlVP-multbdd}
\norm{V_r f}_{\cF\cA} \leq C\norm{f}_{\cF\cA} \qquad\text{ for all $f\in \fourier$\/.}
\end{equation} 

By \cite[Proposition 6.2.5]{Reit_Steg}, $\cF\cA$ contains all
functions in $\fourier$ that
have compact support; in particular, $V_r\in\cF\cA$ for all~$r$\/.
\begin{lem}\label{l:escape}
$\norm{V_r}_{\cF\cA} \to \infty$ as $r\to \infty$\/.
\end{lem}

\begin{proof}
Suppose otherwise; then there would exist a constant $K>0$ and a strictly increasing sequence $r_1<r_2<\ldots$ in $(0,\infty)$\/, such that $\norm{V_{r_n}}_{\cF\cA}\leq K$ for all~$n$\/.
Hence, for any $f\in \fourier$\/, applying \eqref{eq:Segal-condition} yields
\[ \norm{fV_{r_n}}_{\fourier}\leq \norm{fV_{r_n}}_{\cF\cA} \leq C\norm{f}_{\fourier}\norm{V_{r_n}}_{\cF\cA} \leq KC\norm{f}_{\fourier} \qquad\text{ for all $n$\/.} \]
Now, given $f\in C_c(\widehat{\Real^d})$\/, choose $n$ sufficiently large that $\supp(f)\subseteq [-r_n,r_n]^d$\/. Since $fV_{r_n} = f$\/, the preceding inequalities imply that
\begin{equation}\label{eq:Grapelli}
 \norm{f}_{\fourier} \leq \norm{f}_{\cF\cA} \leq KC \norm{f}_{\fourier}\,.
\end{equation} 
Since $C_c(\widehat{\Real^d})$ is dense in $\cF\cA$ (see~\cite[Proposition~6.2.8]{Reit_Steg}), it follows from \eqref{eq:Grapelli} that the norms $\norm{\cdot}_{\cF\cA}$ and $\norm{\cdot}_{\fourier}$ are equivalent, and hence that the norms $\norm{\cdot}_{\cA}$ and $\norm{\cdot}_{L^1(\Real^d)}$ are equivalent. Since $\cA$ is dense in $L^1(\Real^d)$ this implies that $\cA=L^1(\Real^d)$, which is a contradiction.
\end{proof}

For each $n$\/, let $E_n = \cF^{-1}(V_{2^n})$\/. We know by \eqref{eq:dlVP-multbdd} that the sequence $(V_{2^n})$ is bounded in the $\mnorm{\cdot}$-norm, and Lemma~\ref{l:escape} tells us that it is unbounded in $\cF\cA$\/; so pulling back by $\cF$ we see that $(E_n)$ is unbounded but multiplier-bounded in~$\cA$. It remains only to note that since $\cF(E_{n+1})$ is equal to~$1$ on the cube $[-2^{n+1},2^{n+1}]^d$\/, while $\cF(E_n)$ is supported inside that cube, we have $E_nE_{n+1}=E_n$ for all $n$\/. Then, invoking Lemma~\ref{l:absorbing-unbounded} and Theorem~\ref{t:better-nonAA}, we conclude that $\cA$ is not \appam.
\end{proof}

\begin{rems}\
\begin{YCnum}
\item Unlike the examples of the previous section, the sequence $(E_n)_{n\geq 1}$ that was constructed above does not consist of idempotents.
 This illustrates how it would have been too restrictive to only consider orthogonal idempotents in the definition of a SUM configuration.
\item Essentially the same arguments show that proper Segal subalgebras of the convolution algebra $L^1({\mathbb T}^d)$\/, $1\leq d < \infty$\/, are not \appam. We have omitted the details; the result had also been observed by R.~J. Loy, in response to a preliminary draft of this manuscript (private communication).
\item It seems plausible that a similar argument can be carried out for (abstract Segal subalgebras of) other Fourier algebras. The difficulty would lie in obtaining suitable analogues of the \dlVP\ kernel, in order to allow an appeal to Lemma~\ref{l:absorbing-unbounded}.
\end{YCnum}
\end{rems}

\subsection*{Little Lipschitz algebras on compact metric spaces}
Let $(X,d)$ be a compact, infinite metric space, and fix $\al\in(0,1]$\/.
If $f\in C(X)$ let $\Delta f(x,y) = d(x,y)^{-\al}[f(x)-f(y)]$ for $x\neq y$\/, and put
\[ \snorm{f} = \sup_{x\neq y} \abs{\Delta f(x,y)}\,. \]

We briefly recall the definition of the Lipschitz algebras associated to $(X,d)$ and $\al$\/: further details and proofs can be found in~\cite{Sher_lip64}. Recall that $\Lip_\al(X,d)$ is the space of all $f\in C(X)$ such that $\snorm{f}<\infty$\/, equipped with the norm $\norm{f}_\al = \norm{f}_\infty + \snorm{f}$\/; and $\lip_\al(X,d)$ is the closed subspace of $\Lip_\al(X,d)$ consisting of those $f\in \Lip_\al(X,d)$ such that
\[ \Delta f(x,y) \to  0 \text{ as $d(x,y)\to 0$\/.} \]
We equip $\lip_\al(X,d)$ with the norm inherited from $\Lip_\al(X,d)$; with respect to pointwise product of functions, $\lip_\al(X,d)$ is then a unital, commutative Banach algebra, which by results of \cite{BCD_lip} is known not to be amenable.

\begin{thm}\label{t:lip-not-aa}
$\lip_\al(X,d)$ is not \appam.
\end{thm}

Our proof is a slightly indirect application of the previous techniques. Since $X$ is infinite and compact, it contains at least one non-isolated point, $z_0$ say. Let
\[ \cA_0=\{f\in \lip_\al(X,d): f(z_0)=0\}\,;\]
then $\cA_0$ is a maximal ideal in $\lip_\al(X,d)$\/, and the latter may be realized as the unitization of the former.
 Our aim is to show that $\cA_0$ is not \appam, by constructing a sequence $(E_n)\subset\cA_0$ which satisfies the conditions of Lemma~\ref{l:absorbing-unbounded}.

Since $z_0$ is a non-isolated point in $X$, there exists a sequence $(y_n)\subset X\setminus\{z_0\}$ such that $d(y_n,z_0)\to 0$\/. Furthermore, by passing to a subsequence if necessary, we may assume that $d(y_{n+1},z_0) < d(y_n,z_0)/2$ for all~$n$\/.

For ease of notation, we put $\delta_n\defeq d(y_n,z_0)$ and let $V_n\defeq \{ x\in X \st d(x,z_0)\leq \delta_n/2$\}.
Also, given a compact subset $S\subseteq X$ and $x\in X$, we write $d(x,S)$ for the distance from $x$ to~$S$\/, i.e.
\[ d(x,S)\defeq \inf\{ d(x,y)\st y\in S \}\,.\]
(The infimum is attained, by compactness.)

Define $E_n\in C(X)$ by
\begin{equation}\label{eq:bump-function}
E_n(x) \defeq  \min\{ 1, 2\delta_n^{-1} d(x,V_n)\}\,.
\end{equation}


Clearly $0\leq E_n(x)\leq 1$ for all $x$\/, while $E_n$ is zero on~$V_n$\/. Moreover, if $E_n(x)<1$ then there exists $u\in V_n$ with $d(x,u)<\delta_n/2$\/, so that $d(x,z_0)< \delta_n$ by the triangle inequality. Therefore,
\begin{equation}
\text{ if $x\in X$ and $d(x,z_0)\geq \delta_n$ then $E_n(x)=1$\/.}
\end{equation} 
In particular, the hypothesis that $\delta_{n+1}< \delta_n/2$ implies that $E_{n+1}(x)=1$ whenever $x$ lies outside $V_n$\/, and it follows that
\begin{equation}
E_n(x)E_{n+1}(x) = E_n(x) \qquad\text{ for all $x\in X$\/.}
\end{equation}

\begin{lem}
Let $x,y\in X$ with $x\neq y$\/.
Then
\begin{equation}\label{eq:keybound}
\abs{\Delta E_n(x,y)} \leq  \min\{ 2\delta_n^{-1}d(x,y)^{1-\al}, d(x,y)^{-\al} \} \leq 2^\al\delta_n^{-\al}\,. 
\end{equation}
\end{lem}

\begin{proof}
Since $\Delta E_n$ is an antisymmetric function, it suffices to consider the case $d(x,V_n)\geq d(y,V_n)$\/.
Then, if $d(y,V_n)\geq \delta_n/2$ we have $E_n(x)=E_n(y)=1$\/, so that $\Delta E_n(x,y)=0$\/.

We thus restrict attention to the case where $d(y,V_n) < \delta_n/2$\/. For such~$y$\/, we have $E_n(y)=2\delta_n^{-1}d(y,V_n)$\/. Then, since $E_n(x)\leq 2\delta_n^{-1}d(x,V_n)$ for all $x$\/, and since $d(x,V_n)\leq d(y,V_n)+d(x,y)$ by the triangle inequality, we have
\begin{subequations}
\begin{equation}\label{eq:local}
0 \leq \Delta E_n(x,y) \leq \frac{2\delta_n^{-1}d(x,V_n)-2\delta_n^{-1}d(y,V_n)}{d(x,y)^\al}
 \leq 2\delta_n^{-1} d(x,y)^{1-\al}\,.
\end{equation}
On the other hand, we have a trivial upper bound
\begin{equation}\label{eq:spread}
0 \leq \Delta E_n(x,y) \leq \frac{1-0}{d(x,y)^\al} = d(x,y)^{-\al}\,.
\end{equation}
\end{subequations}
Combining \eqref{eq:local} and \eqref{eq:spread} gives the first inequality in \eqref{eq:keybound}. The second one follows by observing that, since $0<\al\leq 1$\/, the function $t\mapsto \min\{2\delta_n^{-1} t^{1-\al}, t^{-\al}\}$ attains its maximum at $t=\delta_n/2$\/.
\end{proof}

\begin{prop}\label{p:key-sequence-for-lip}
Let $n\in\Nat$\/. Then $E_n\in\cA_0$\/, and:
\begin{YCnum}
\item $\delta_n^{-\al} \leq \snorm{E_n} \leq 2^\al\delta_n^{-\al}$\/;
\item for all $f\in\cA_0$\/, we have $\norm{E_nf}_\infty\leq \norm{f}_\infty$ and $\snorm{E_nf}\leq 3\snorm{f}$\/.
\end{YCnum}
\end{prop}

\begin{proof}
The inequalities in \eqref{eq:keybound} show that $E_n\in \lip_\al(X,d)$\/, and by construction $E_n(z_0)=0$\/. Thus $E_n\in\cA_0$\/. The upper bound in (i) is also immediate from \eqref{eq:keybound}, while the lower bound follows from the estimate
\[ \snorm{E_n}\geq \Delta E_n(y_n,z_0) = \frac{2\delta_n^{-1}d(y_n,V_n)}{d(y_n,z_0)^\alpha} \geq 2\delta_n^{-1}\frac{d(y_n,z_0)- \delta_n/2}{d(y_n,z_0)^\alpha} = \delta_n^{-\alpha} \,.\] 

The first estimate in (ii) is trivial. For the second, we argue as follows. If both $d(x,V_n)$ and $d(y,V_n)$ are $\geq \delta_n/2$\/, then it follows from the definition of $E_n$\/, that
\[  \abs{\Delta (E_n f)(x,y)} = \abs{\Delta f(x,y)}\leq \snorm{f}\,.\]
If $d(x,V_n) < \delta_n/2$\/, then $d(x,z_0)\leq\delta_n$\/, and so since $f(z_0)=0$ we have
\[ \abs{f(x)} = d(x,z_0)^\al \abs{\Delta f (x,z_0)} \leq {\delta_n}^\al \snorm{f} \,.\]
Hence using part (i) yields
\[ \begin{aligned}
\abs{\Delta (E_nf)(x,y)}
 & = \abs{ [\Delta E_n](x,y)\cdot f(x) + E_n(y)\cdot[\Delta f](x,y) }\\
 & \leq \snorm{E_n} \cdot {\delta_n}^\al \snorm{f} + 1\cdot \snorm{f} \\
 & \leq 3 \snorm{f}\,.
\end{aligned} \]
By symmetry, the same argument with $x$ and $y$ interchanged shows that if $d(y,V_n) < \delta_n/2$ then $\abs{\Delta(E_nf)(x,y)}\leq 3\snorm{f}$\/. Thus $\snorm{E_n f}\leq 3\snorm{f}$, as claimed.
\end{proof}

\begin{proof}[Proof of Theorem~\ref{t:lip-not-aa}]
By Proposition~\ref{p:key-sequence-for-lip}, we can apply Lemma~\ref{l:absorbing-unbounded} and thence Theorem~\ref{t:better-nonAA} to the Banach algebra $\cA_0$, to conclude that it is not \appam. Since there is an isomorphism of Banach algebras $\fu{\cA_0}\iso\lip_\al(X,d)$\/, it follows from \cite[Proposition~2.4]{GhL_genam1} that $\lip_\al(X,d)$ is not~\appam.
\end{proof}

\end{section}

\begin{section}{Approximate amenability of tensor products}
It was claimed in \cite[Proposition 2.3]{GhL_genam1} that if $\cA$ is \appam\ and has a bounded approximate identity, and $\cB$ is amenable, then $\cA\ptp \cB$ is \appam\/. However, the argument there is incomplete; the correct partial result is as follows.

\begin{prop}\label{p:amended_GhL}
Let $\cA$ and $\cB$ be Banach algebras, where $\cA$ is \appam\ and has a bounded approximate identity, and $\cB$ is amenable. Let $X$ be a Banach $\cA\ptp\cB$-bimodule.
Then for every bounded derivation $D:\cA\ptp\cB \to X^*$\/, there exists a net $(\xi_\al)\subset X^*$ such that $D(a\tp b) = \lim_\al (a\tp b)\cdot \xi_\al -\xi_\al\cdot (a\tp b)$ for all $a\in\cA$ and all $b\in\cB$\/.

Moreover, if $\cA$ is furthermore assumed to be \emph{boundedly} \appam, then
\begin{YCnum}
\item there exists $C>0$ such that the net $(\xi_\al)$ may be chosen to satisfy
\begin{equation}\label{eq:saving_bound}
 \norm{(a\tp b)\cdot\xi_\al -\xi_\al\cdot(a\tp b)}\leq C\norm{a}\norm{b} \qquad\text{ for all $a\in\cA$\/, $b\in\cB$\/;}
\end{equation} 
\item $\cA\ptp\cB$ is boundedly \appam.
\end{YCnum}
\end{prop}

\begin{proof}
For the first part, we follow the proof of \cite[Proposition 5.4]{BEJ_CIBA}: note that in this argument it is important to first use amenability of $\cB$ and only then the \appamy\ of $\cA$\/, rather than the other way round. If we assume furthermore that $\cA$ is boundedly \appam, then inspection of the argument shows that the net $(\xi_\al)$ thus produced will satisfy the upper bound claimed in~\eqref{eq:saving_bound}. Finally, using \eqref{eq:saving_bound} and the definition of the projective tensor product, it is easily checked that $\norm{w\cdot\xi_\al-\xi_\al\cdot w} \leq C\norm{w}$ and $D(w)=\lim_\al w\cdot\xi_\al-\xi_\al\cdot w$ for all $w\in \cA\ptp\cB$\/; thus, every bounded derivation from $\cA\ptp\cB$ to a neo-unital bimodule is the strong-operator limit of a bounded net of inner derivations. By the reasoning of \cite[Proposition~2.5]{GhL_genam1}, this suffices to show that $\cA\ptp\cB$ is boundedly \appam.
\end{proof}

We do not know if the tensor product of two (boundedly) \appam\ Banach algebras is itself (boundedly) \appam. The following results allow us to make some progress.

\begin{lem}\label{l:tensor-with-findim}
Let $\cA$ be \appam\ and let $\cB$ be a finite-dimensional, amenable algebra. Then $\cA\ptp\cB$ is \appam.
\end{lem}
\begin{proof}
It suffices to consider derivations $D:\cA\ptp\cB\to X^*$\/, where $X$ is a neo-unital $\cA\ptp\cB$ bimodule $X$\/. By Proposition~\ref{p:amended_GhL}, there exists a net $(\xi_\al)\subset X^*$ such that, for every $a\in\cA$ and $b\in\cB$\/,
$D(a\tp b) = \lim_\al (a\tp b)\cdot \xi_\al - \xi_\al (a\tp b)$\/.
But since $\cB$ is finite-dimensional, every element of $\cA\ptp\cB$ can be written as a \emph{finite} linear combination of elementary tensors, so by linearity we have $D(w)=\lim_\al w\cdot\xi_\al-\xi_\al\cdot w$ for every $w\in \cA\ptp\cB$\/, and the proof is complete.
\end{proof}

\begin{thm}\label{t:tensor-product}
Let $\cA$ and $\cB$ be Banach algebras,
with \emph{central} bounded approximate identities $(u_i)_{i\in I}$ and $(v_j)_{j\in J}$ respectively. Suppose that there exists a directed family $(\cB_\gamma)_{\gamma\in\Gamma}$ of finite-dimensional, amenable subalgebras of $\cB$\/,
that satisfy the following conditions:
\begin{YCnum}
\item $\bigcup_\gamma \cB_\gamma$ is dense in $\cB$\/;
\item for each $j$\/, there exists $\gamma(j)$ such that $v_j\cB = \cB v_j\subseteq \cB_{\gamma(j)}$\/.
\end{YCnum}
Then if $\cA$ is \appam, so is $\cA\ptp\cB$\/.
\end{thm}

\begin{proof}[Proof of theorem]
Equip $I\times J$ with the product ordering: then the net $(u_i\tp v_j)_{I\times J}$ is a central, bounded approximate identity for $\cA\ptp\cB$\/, and
$(u_i\tp v_j)\cB = \cB(u_i\tp v_j) \subseteq \bigcup_{\gamma} \cA\ptp\cB_\gamma$ for every $(i,j)\in I\times J$\/.
Also, by Lemma \ref{l:tensor-with-findim}, for each $\gamma$ the algebra $\cA\ptp\cB_\gamma$ is approximately amenable.
Therefore, by \cite[Theorem 2.3]{GhSt_AG}, $\cA\ptp\cB$ is pseudo-amenable; since it has a bounded approximate identity, it is \appam, by \cite[Proposition~3.2]{GhZ_pseudo}.
\end{proof}

\begin{eg} In Theorem \ref{t:tensor-product} we can take $\cB$ to be any of the following examples:
\begin{YCnum}
\item $\cB= C_r^*(G)$ where $G$ is a compact group. In this case the dual object $\widehat{G}$ is discrete, and we may identify $\cB$ with the $c_0$-sum $\bigoplus_{\sigma\in\widehat{G}} {\mathbb M}_{d(\sigma)}(\Cplx)$\/, where $d(\sigma)$ is the (finite) dimension of the irreducible representation $\sigma$\/.
Take $\Gamma$ to be the set of finite subsets of $\widehat{G}$\/, partially ordered by inclusion. For each $F\in\Gamma$ let $\cB_F$ be the $c_0$-sum  $\bigoplus_{\sigma\in F} {\mathbb M}_{d(\sigma)}(\Cplx)$\/, and let $v_F$ be the identity element of the finite-dimensional algebra $\cB_F$\/.  Then $(v_F)_{F\in\Gamma}$ is a central bounded approximate identity in $\cB$\/, and it is easily verified that the conditions of Theorem~\ref{t:tensor-product} are satisfied.
\item $\cB=\lp{1}(\Nmin)$\/, the convolution algebra of the semilattice $\Nmin$\/. Denoting the canonical unit basis of $\lp{1}(\Nmin)$ by $(e_j)_{j\in\Nat}$\/, let $\cB_n$ be the linear span of $\{e_1,\ldots, e_n\}$\/. Then $\cB_n$ is a subalgebra (in fact, an ideal) in~$\cA$ and
 is known to be isomorphic to $\Cplx^n$ with pointwise multiplication (so in particular is amenable). Moreover, for each~$n$ we have $e_n\cB = \cB_n = \cB e_n$\/, and $(e_n)_{n\geq 1}$ is a bounded approximate identity for $\cB$\/. Hence, once again, the conditions of Theorem~\ref{t:tensor-product} are satisfied.
\end{YCnum}
\end{eg}
\end{section}

\begin{section}{Commutative algebras whose underlying Banach space is reflexive}
Any amenable Banach algebra whose underlying Banach space is Hilbertian (i.e.~isomorphic to a Hilbert space) must be finite-dimensional: this is well-known and follows from \cite[Theorem~2.2]{BEJ_WCH92}. Various authors have generalized this result to other reflexive Banach spaces, but only under certain extra hypotheses on the ideal structure of the Banach algebra in question.

In this section we make some partial progress towards corresponding results for bounded \appamy. Specifically, we obtain a complete answer in the case of \emph{commutative} Banach algebras. Our main tool is the following result.

\begin{prop}\label{p:TERRY}
Let $\cA$ be a unital, commutative Banach algebra. Suppose that every maximal ideal in $\cA$ has an identity element. Then $\cA$ is finite-dimensional and isomorphic to $\Cplx^N$ for some~$N$\/.
\end{prop}

\begin{proof}
Let $\Phi$ be the character space of $\cA$, equipped with the Gelfand topology. Our first step is to show that $\Phi$ is finite; since $\Phi$ is compact, it will suffice to show that each point of $\Phi$ is isolated (i.e.~open and closed).

Let $\varphi\in\Phi$\/. By assumption, $\ker\varphi$ has an identity element, $u_\varphi$ say; a small calculation yields
\begin{equation}\label{eq:TERRY}
(1-u_\varphi) a = \varphi(a) (1-u_\varphi) \qquad\text{ for all $a\in \cA$\/.}
\end{equation}
Now, if $\chi\in\Phi\setminus\{\varphi\}$ then there exists $b\in \cA$ such that $\chi(b)\neq\varphi(b)$\/. Hence, by \eqref{eq:TERRY},
\[ (1-\chi(u_\varphi))\chi(b) = \chi((1-u_\varphi) b) = (1-\chi(u_\varphi)) \varphi(b) \]
so that $\chi(u_\varphi)=1$\/. Hence the Gelfand transform of $u_\varphi$ is zero at $\varphi$ and $1$ at all other points of $\Phi$\/, which implies that $\varphi$ is open and closed in $\Phi$ as required.

Therefore, $\Phi$ is a finite set, consisting of characters $\varphi_1,\ldots,\varphi_N$\/, say. Using the elements $1-u_{\varphi_j}$\/, for $1\leq j \leq N$\/,it is clear that the Gelfand transform is a surjection onto $C(\Phi)\iso\Cplx^N$\/.
To finish, note that the Jacobson radical of $\cA$\/, which we shall denote by $\Rad(\cA)$\/, is equal to $\bigcap_{j=1}^N \ker\varphi_j$\/. Put $u \defeq u_{\varphi_1} \dotsb u_{\varphi_N}\in\Rad(\cA)$\/; then $ux=x$ for all $x\in\Rad(\cA)$. This is only possible if $\Rad(A)=0$\/, and so the Gelfand transform is injective, which completes the proof. 
\end{proof}

Recall that a Banach algebra $\cA$ is a \dt{dual Banach algebra}, if $\cA = (\cA_*)^*$ for some Banach space $\cA_*$\/, and the algebra multiplication is separately weak* continuous. We denote the canonical embedding of $\cA_*$ into $\cA^*$ by $\imath$\/. 
\begin{lem}\label{l:id-in-DBA}
Let $\cA$ be a dual Banach algebra and be boundedly \appam . Then $\cA$ has an identity.
\end{lem}

\begin{proof}
From \cite[Theorem 2.5]{CGZ_JFA09} there are nets $(F_i), (G_i) \subseteq \cA^{**}$\ such that 
$aF_i\rightarrow a$ and $G_ia \rightarrow a$,\ for every $a\in \cA$, and $(F_i)$ and $(G_i)$ are multiplier bounded.  Since $\cA=(\cA_*)^*$ is a dual Banach algebra, $\imath^*$
is an $\cA$ bimodule morphism from $\cA^{**}$ onto $\cA$, and so $(\imath^*(F_i))$\ and $(\imath^*(G_i))$ are 
multiplier-bounded right and left approximate identities for~$\cA$. Therefore, by \cite[Theorem 3.3]{CGZ_JFA09} $\cA$ has a bounded approximate identity; taking a $\wstar$ cluster point yields an identity.
\end{proof}

\begin{thm}
Let $\cA$ be a commutative, boundedly \appam\ Banach algebra. If the underlying Banach space of $\cA$ is reflexive, then $\cA$ is finite-dimensional.
\end{thm}

\begin{proof}
By the previous lemma,
$\cA$ has an identity element~$\id[\cA]$\/. By renorming $\cA$ if necessary, we may assume that $\norm{\id[\cA]}=1$\/, i.e.~that $\cA$ is unital.

Now let $M$ be a maximal ideal in $\cA$\/. Clearly $M$ is reflexive; moreover, since $\fu{M}\iso \cA$ and $\cA$ is boundedly \appam, $M$ is also boundedly \appam\ (see \cite[Proposition~2.4]{CGZ_JFA09}), and hence (by Lemma~\ref{l:id-in-DBA})
it has an identity element. Applying Proposition~\ref{p:TERRY} completes the proof.
\end{proof}
\end{section}

\begin{section}{Bounded approximate amenability of ${\rm A}(G)^{**}$ and ${\rm A}(G)$}
For $G$ a locally compact group, $L^1(G)$ denotes the group algebra; $\FA(G)$ and $\FS(G)$ denote, respectively, the \dt{Fourier} and \dt{Fourier-Stieltjes} algebra of~$G$\/, as defined by Eymard~\cite{Eym_BSMF64}.
Amenability and \appamy\ of $L^1(G)^{**}$ were studied in \cite{GhLW_PAMS96} and \cite{GhL_genam1} respectively, while amenability of $\FA(G)^{**}$ was characterized in \cite{Gran_JAust97} and \cite{LauLoy_JFA97}.
Here we characterize bounded \appamy\ of~$\FA(G)^{**}$\/.

Throughout this section, the second dual algebras are equipped with the first Arens product, which we denote by~$\boxprod$. For further details on Arens products, see \cite[\S3.3]{Dal_BAAC} or \cite[\S1.4]{Palmer1}.

\begin{prop}\label{p:technical-for-A**}
Let $\cA$ be a commutative Banach algebra. If $\cA^{**}$ is boundedly \appam, then $\cA$ has a bounded approximate identity. Moreover, $\cA^{**}$ has an identity.
\end{prop}

\begin{proof}
From Theorem~2.6 of \cite{CGZ_JFA09} there are nets 
$(M_i)\subseteq  (\cA^{**} \widehat{\otimes}\cA^{**})^{**}$ , 
$(F_i)\subseteq \cA^{****}$, $(G_i)\subseteq \cA^{****}$
and $K>0$\/, such that 
\begin{equation}\label{eq:originally0.1}
\begin{aligned} 
& \norm{ u\cdot M_i - M_i\cdot u - u\otimes G_i +F_i \otimes 
{u} } \leq K \norm{u}, \\  &u\cdot F_i\rightarrow u,\ \  
 \norm{ u\cdot F_i} \leq K \norm{ u},  \\
& G_i\cdot u\rightarrow u,\ \
    \norm{G_i\cdot u} \leq K \norm{ u},\ \  ( u\in \cA^{**}).\end{aligned}
\end {equation} 
Let $\imath$ be the canonical embedding of $\cA^*$ in $\cA^{***}$\/. 
Then $P=\imath^*$ is an $\cA$- module morphism from $\cA^{****}$ onto $\cA^{**}$. Now, by using 
weak*-continuity of the first Arens product in the first variable, 
we can show that $P$ is a right $\cA^{**}$-morphism. In fact, take $F\in \cA^{****}$ and $m\in \cA^{**}$. Let $(f_i)$ be a net
of elements of $\cA^{**}$ converging weak* to $F$. Then
\[ P(F\cdot m)  = \wstar\lim_i P(f_im) = \wstar\lim_i f_i m = \wstar\lim_i P(f_i)m = P(F)m\,. \]
Hence by applying $P$ to the third equation in \eqref{eq:originally0.1}
 we see that $\cA^{**}$ has a multiplier-bounded left approximate identity, $(f_\beta)$, say.
From \eqref{eq:originally0.1} again, for $u=f_\beta$ we have
\begin{equation}\label{eq:originally0.2}
  \norm{(f_\beta\cdot M_i- M_i f_\beta- f_\beta \otimes G_i + F_i\otimes 
f_\beta)\cdot F_\al} \leq K^2 \norm{f_\beta},\;
 \text{for every $\al$, $\beta$ and $i$\/.}
\end{equation}
By using the triangle inequality, and  left-multiplier-boundedness of the set $(f_\beta)$\/, from \eqref{eq:originally0.2} we have             
\begin{equation}
\begin{aligned}
   \norm{ f_\beta} \norm{ G_i \boxprod F_\al }
 & \leq K^2 \norm{f_\beta} +
 \norm{f_\beta\cdot M_i\cdot F_\al}+ \norm{M_i\cdot(f_\beta F_\al)}+ \norm{F_i} \norm{f_\beta F_\al}     \\
 & \leq 
K^2 \norm{f_\beta}+K \norm{M_i\cdot F_\al}+K \norm{M_i} \norm{F_\al}+K \norm{F_i}
 \norm{F_\al},     \\
 & \text{for all $\al$\/, $\beta$\/ and $i$\/.}\end{aligned}
\end{equation} 
Hence,
\begin{equation}\label{eq:originally0.4}
 \norm{G_i\boxprod F_\al}  \leq 
K^2+\frac{K}{ \norm{f_\beta}}( \norm{M_i\cdot F_\al}+ \norm{M_i} \norm{F_\al}
+ \norm{F_i} \norm{F_\al}),
 \text{for all $\al$\/, $\beta$ and $i$\/.}
\end{equation}

Now two cases may occur.

\paragraph{Case 1:} $(f_\beta)$ is unbounded.
In this case, from \eqref{eq:originally0.4} we have
\[  \norm{G_i\boxprod F_\al} \leq K^2, \quad\text{for all $i$ and $\al$\/.} \]
 
For every $f\in\cA$\/, noting that $f\cdot G_i=G_i\cdot f$ as $\cA$ is commutative, we have
\break
$f\cdot G_i\boxprod F_\al\rightarrow f$ as $\al$ and $i$ tend 
to infinity. Hence $f\cdot P(G_i\boxprod F_\al)\rightarrow f$\/.
Taking $E$ to be any weak*-cluster point of $(P(G_i\boxprod F_\al))_{i,\al}$\/, it follows that $E$\ is a right identity for $\cA^{**}$\/, and hence $\cA$ has a bounded approximate identity.

\paragraph{\bf Case 2:} the net $(f_\beta)$ is bounded. Then we may argue as in case~1, to conclude that $\cA$ has a bounded approximate identity.

\medskip For the final part: let $(e_j)$ be a bounded approximate identity for $\cA$, and let $E$ be a 
weak*-cluster point of $(\hat{e_j})$. Then $E$ is a right identity for $\cA^{**}$\/. On the other hand, since $\cA^{**}$ is boundedly \appam, it is in particular \appam. By \cite[Lemma~2.2]{GhL_genam1}, it therefore has a left approximate identity, say $(n_j)$, and we obtain
\begin{equation}
E \boxprod F = \lim_j n_j \boxprod E \boxprod F = \lim_j n_j \boxprod F = F\,,\quad\text{ for every $F\in\cA^{**}$\/.}
\end{equation}
Thus $E$ is  also a left identity for $\cA^{**}$\/. 
\end{proof}

The proof of the following theorem uses, in part, an idea from the proof of \cite[Theorem 3.3]{GhL_genam1},
but also some properties specific to the Fourier algebra. Proposition~\ref{p:technical-for-A**} plays an important role.

\begin{thm}\label{t:BAA-AG**}
Let $G$ be a locally compact group. Then $\FA(G)^{**}$ is boundedly \appam\ if and only if $G$ is finite.
\end{thm}

\begin{proof}
Sufficiency is trivial. Suppose, then, that $\FA(G)^{**}$ is boundedly \appam.

Let $n$ be a topological invariant mean on $\VN(G) = \FA(G)^*$\/, that is, a state which satisfies $f\cdot n = \langle f,I\rangle n$ for all $f\in \FA(G)$, where $I\in\VN(G)$ is the identity operator. Then 
$r\boxprod n = n$ for every state $r\in \VN(G)^*$\/, by $\wstar$-continuity of $\boxprod$ in the first variable. Hence
$J_n = n\boxprod \FA(G)^{**}$ is a closed, complemented, two-sided ideal in  $\FA(G)^{**}$ and so, by~\cite[Corollary 2.4]{GhL_genam1}, $J_n$ has a right approximate identity, $(p_i)$ say.

If $m$ is \emph{any} topological invariant mean on $\VN(G)$\/, then $m=n\boxprod m\in J_n$. Since $p_i\in J_n$\/, $m\boxprod p_i=p_i$\/, and therefore
\begin{equation}
m = \lim m\boxprod p_i =\lim p_i,
\end{equation}
showing that $\VN(G)$ has a unique topological invariant mean. Hence by \cite[Corollary~4.11]{LauLos_JFA93}, $G$~is discrete. 

On the other hand, since $\FA(G)^{**}$ is boundedly \appam, it follows from the preceding proposition that it has a two-sided identity.
By \cite[Proposition~3.2(b)]{Lau_AG**}, $G$ must be compact.
 Since we have also shown that $G$ is discrete, it follows that $G$ is finite.
\end{proof}

We now examine \dt{essential amenability} (\cite[\S1]{GhL_genam1}) of~$\FA(G)^{**}$\/. It suffices to know the following two results: if $A$ is an essentially amenable Banach algebra and there exists a continuous algebra homo\-morphism of $A$ \emph{onto} another Banach algebra~$B$\/, then by \cite[Proposition~2.2]{GhL_genam1}, $B$ is essentially amenable; and if $A$ is an essentially amenable Banach algebra with an identity element, then $A$ is amenable.

Recall that $\UCBhat{G}$ is defined to be the closure of 
$\FA(G)\cdot \VN(G)$, where $\VN(G)$ is identified with the dual of $\FA(G)$ and the action of $\FA(G)$ on $\VN(G)$ is the canonical action of an algebra on its dual.
The reduced group $\Cstar$-algebra $\Cstar_r(G)$ is contained in $\UCBhat{G}$ (see \cite[Proposition~4.4]{Lau79_UCB} for a proof).
 
For the proof of the theorem which follows, it is convenient to work with the notion of a `left-introverted' subspace of $\VN(G)$ (for the definition, see \cite[\S5]{Lau79_UCB}, for example).
Examples of left-introverted subspaces of $\VN(G)$ are: $\VN(G)$ itself; $\UCBhat{G}$; and the reduced $C^*$-algebra $\Cstar_r(G)$ (\cite[Proposition~5.2]{Lau79_UCB}). The relevance of these subspaces to the present study of $\FA(G)^{**}$ arises from the following lemma, which follows from the discussion on~\cite[p.~177]{LauLoy_JFA97}.

\begin{lem}
Let $Y\subseteq X$ be two left-introverted subspaces of $\VN(G)$\/. Then the duals $Y^*$ and $X^*$ can both be equipped with Arens-type products, with respect to which the canonical quotient maps $\VN(G)^*\to X^*$ and $X^*\to Y^*$ are algebra homomorphisms.
\end{lem}

In particular, there are quotient homomorphisms $\FA(G)^{**} \to \UCBhat{G}^*$ and $\UCBhat{G}^*\to \Cstar_r(G)^*$\/. Moreover, the algebra $\Cstar_r(G)^*$ is naturally embedded as a closed subalgebra of $\FS(G)$ (\cite[Propo\-sition~5.3]{Lau79_UCB}), and is \emph{equal} to $\FS(G)$ if and only if $G$ is amenable.

\begin{thm}
Let $G$ be an amenable group. Then either of the Banach algebras $\FA(G)^{**}$ or $\UCBhat{G}^*$ is essentially amenable if and only if $G$ is finite.
\end{thm}

\begin{proof}
If $G$ is finite, then $\FA(G)^{**}=\UCBhat{G}^*=\FA(G)$ is amenable, and in particular essentially amenable.
Conversely, if either $\FA(G)^{**}$ or $\UCBhat{G}^*$ is essentially amenable, then so is $\Cstar_r(G)^*$\/. Since $G$ is amenable, it follows that $\FS(G)=\Cstar_r(G)^*$ is also essentially amenable, and hence amenable as it has an identity. Corollary~2.4 of \cite{ForRun_amenAG} then implies that $G$ is compact, so that $\FA(G)$ has an identity element. Therefore $\UCBhat{G}= \FA(G)\cdot \VN(G)=\VN(G)$\/, and so $\UCBhat{G}^*=\FA(G)^{**}$ has an identity and is essentially amenable, hence is amenable. 
By the results of \cite{Gran_JAust97} (or 
Theorem~\ref{t:BAA-AG**} 
above), we conclude that $G$ is finite.
\end{proof}

For a general locally compact group $G$\/, there is still no complete characterization of the bounded approximate amenability of $\FA(G)$\/. However, in the special case where $G$ has an open abelian subgroup (in particular, if $G$ is discrete), we can give a complete answer.

\begin{thm} Suppose that $G$ has an open abelian subgroup. Then $\FA(G)$ is boundedly approximately amenable, if and only if $G$ is amenable.\end{thm}

\begin{proof} Suppose that $\FA(G)$ is boundedly \appam.
Then from  \cite[Proposition~3.6]{CGZ_JFA09} it follows that every multiplier-bounded subset of $\FA(G)$ is norm bounded.
This is equivalent to saying that the multiplier norm on $\FA(G)$ is equivalent to its  norm. Hence the amenability of $G$ follows from \cite[Theorem 1]{Losert_PAMS84}. The converse follows from a direct application of \cite[Corollary~3.2(ii)]{GhSt_AG}.
\end{proof}
\end{section}

\begin{section}{Concluding remarks}
We have seen that many of the existing `non-amenability' results for Banach algebras arising in abstract harmonic analysis and operator theory admit extensions to the (boundedly) \appam\ setting; however, in some of these cases, the proofs require genuinely new arguments rather than mere extension of the existing ones.

The technique presented in Section~\ref{s:criterion} seems deserving of further exploration. Several of the technical hypotheses can perhaps be weakened, although it is unclear if doing so would give a substantially stronger result. The applications presented in Sections \ref{s:criterion}~and~\ref{s:application} are intended to be illustrative, rather than exhaustive, and we hope they will stimulate further developments.

\subsection*{Acknowledgements}
The authors would like to thank R. Stokke for useful discussions on Fourier and Fourier-Stieltjes algebras, and Y. Zhang for a careful reading of the present work. 

Some of this work was undertaken during a visit of the first author to the University of Manitoba in April 2009. He thanks the University of Manitoba for their support and hospitality.
\end{section}



\vfill

\noindent%
\begin{tabular}{l@{\hspace{20mm}}l}
Y. Choi
	&  F. Ghahramani \\
D\'epartement de math\'ematiques
	& Department of Mathematics, \\
 \text{\hspace{1.0em}} et de statistique,
	& \\
Pavillon Alexandre-Vachon
	& Machray Hall\\
Universit\'e Laval & University of Manitoba\\
Qu\'ebec, QC
	& Winnipeg, MB \\
Canada, G1V 0A6 & Canada, R3T 2N2 \\
	& \\
{\bf Email: \tt y.choi.97@cantab.net}
	&{\bf Email: \tt fereidou@cc.umanitoba.ca} 
\end{tabular}

\end{document}